\documentclass{article}
\usepackage{graphicx} 
\usepackage{amsmath,amsfonts, amsthm, amssymb, xcolor, enumitem}
\usepackage{tabularx, longtable, array, verbatim}
\usepackage{mathtools}
\usepackage[maxbibnames=8,url=true,style=alphabetic,doi=false]{biblatex}
\usepackage{geometry}
\usepackage{subcaption}

\usepackage[vlined,ruled,algo2e]{algorithm2e}
\SetAlCapHSkip{0em}
\SetCustomAlgoRuledWidth{\textwidth}

\bibliography{survey}

\newtheorem{theorem}{Theorem}

\makeatletter
\let\c@algocf\c@theorem

\makeatother

\newcommand{\Z}{\mathbb{Z}}

\newcommand{\T}{\mathbb{T}}

\newcommand{\C}{\mathbb{C}}
\newcommand{\R}{\mathbb{R}}

\title{Nonlinear Fourier series}
\author{
  Michel Alexis \\ School of Math. and Stat. Sciences\\ Clemson University \\ \texttt{micalexis.math@gmail.com}
  \and
  Christoph Thiele \\ Mathematical Institute\\ Bonn University \\ \texttt{thiele@math.uni-bonn.de}
}

\date{}

\begin{document}

\maketitle

\begin{abstract}
    We survey aspects of nonlinear Fourier series as presented in the 2026 Brouwer lecture 
    by the second author. The fast nonlinear Fourier expansion mentioned in the lecture is further 
    elaborated in joint work with the first author.
\end{abstract}
\section{Fourier series}

About two hundred years ago, Joseph Fourier
advertised in his work \cite{MR2856180} a series expansion for functions, now named after him, as a tool to solve partial differential equations. 
He claimed that Fourier series  can represent general functions on an interval, sparking decades of research into the nature of general functions and their representations.
Over time, many surprising applications of Fourier series appeared. Dirichlet led this way with his  breakthrough \cite{dirichlet1837beweis} in analytic number theory, establishing the expected density of primes in arithmetic progressions.

A finite Fourier series represents a complex-valued function $v$ in a variable $\theta\in [0,2\pi)$ as
\begin{equation}\label{fourier}v(\theta)
    =\sum_{n=-N}^N
    f_n e^{i n\theta} 
    =\sum_{n=-N}^N
    f_n z^n\, .
 \end{equation} 
By Euler's formula, the exponential $e^{i n\theta}$ can be written as a linear combination of sine and cosine functions, 
elucidating the classical interpretation of the Fourier series as a superposition of pure frequencies. Writing $z$ for the complex number $e^{i\theta}$ in the unit circle $\T$, the Fourier series becomes a Laurent polynomial in $z$. Therefore, finite Fourier series are also called trigonometric polynomials.

Trigonometric polynomials are insufficient to represent very general functions, in particular indicator functions of non-trivial subsets of $\T$. This suggests taking
a limit of \eqref{fourier} as $N$ tends to infinity, as Fourier proposed.
This limit raises the questions of convergence to the 
limit and of the class of functions that can be represented with a particular notion of convergence.

If the coefficients $f_n$ are absolutely summable, 
\begin{equation}
    \sum_{n\in \Z} |f_n|<\infty\, ,
\end{equation}
then the Fourier series \eqref{fourier} converges uniformly to a continuous function as $N$ tends to infinity. One still fails to obtain any non-trivial indicator function in this way. Not even every continuous function can be represented as an absolutely summable Fourier series, but only functions in what is called the Wiener space. 

If the coefficients $f_n$ are square summable,
\begin{equation}\label{l2}\sum_{n\in \Z} |f_n|^2<\infty\ ,
\end{equation}
 then the series \eqref{fourier} does not necessarily have a limit at every point 
$\theta$ as $N$ tends to infinity. However, the series converges in the Hilbert space sense,
\begin{equation}
    \lim_{N\to \infty} \int_0^{2\pi} \left|v(\theta)- \sum_{n=-N}^N f_ne^{in\theta}\right|^2\, d\theta=0\, .
\end{equation}
This allows one to represent a measurable and square integrable function $v$, and in particular
the indicator function of a measurable set. The theory of Hilbert spaces was partially motivated by the desire to clarify the problem of convergence of Fourier series. 

\begin{figure}
    \includegraphics[width=\textwidth]{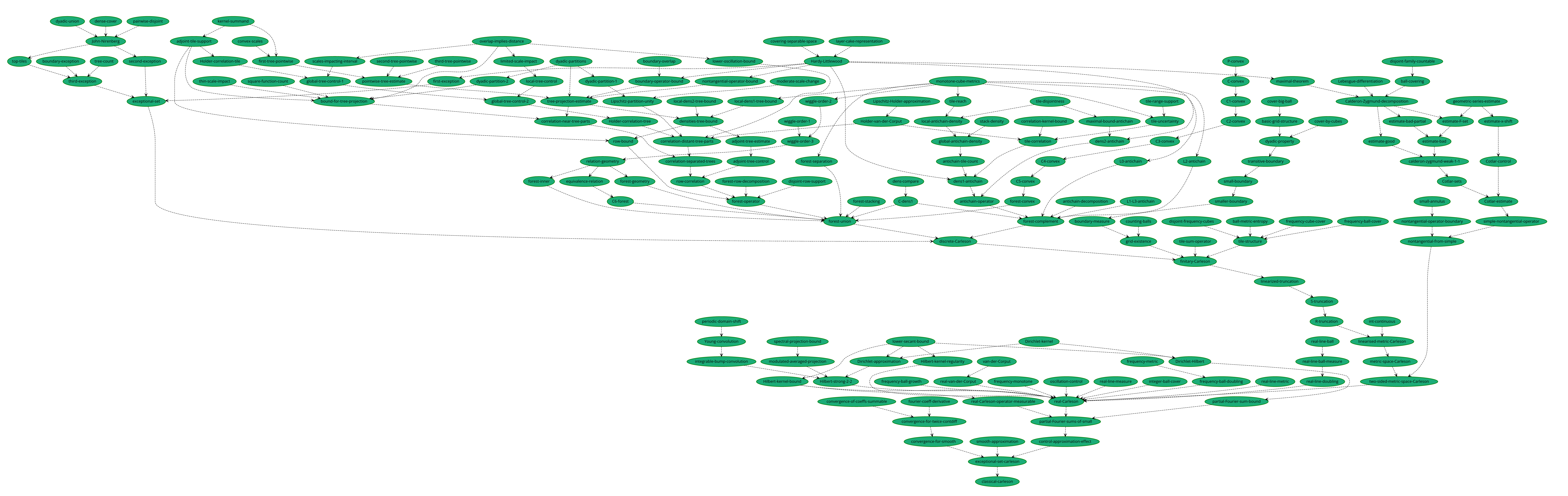}
    \caption{The Lean Carleson project}
    \label{lc}
\end{figure}

In seminal work, Carleson \cite{MR199631} showed that for square summable coefficients
the series
\eqref{fourier} converges  pointwise for almost every $\theta$ as $N$ tends to infinity.
This celebrated theorem was chosen for a formalization project \cite{B26}. Formalizing a theorem means writing it and its proof in code
that allows the proof to be checked for correctness by a computer. Carleson's theorem
was split into about 200 pieces, see Figure \ref{lc}, which were
coded by 12 contributors within a year.
At the time of the project between 2024 and 2025, experts were surprised by the speed of the distributed effort to formalize such a difficult theorem in analysis. Less than two years later, similar formalization tasks can be achieved automatically by artificial intelligence in a few days.
The computer language Lean is becoming a standard for such formalizations. Lean is a successor of Automath, the first 
successful computer verification system developed by N.G.\,(Dick) de Bruijn
\cite{MR1429395}. De Bruijn was inspired by the intuitionistic mathematics of L.E.J.\,(Bertus) Brouwer. 

Fourier introduced his series as part of a solution method to the heat equation.
We present his idea in the case of the related Airy equation 
\begin{equation}\label{airy}
\partial_t u(\theta,t)= \partial_\theta^3 u(\theta,t)    
\end{equation}
for a function in the variable 
 $\theta$ and a time variable $t$. The typical question for this equation
concerns its evolution in time. We are given an initial value of $u$ at time $t=0$ in the form
\begin{equation}\label{initial}
u(\theta,0)= v(\theta)    
\end{equation}
and seek to extend $u$ to all times $t\ge 0$ satisfying \eqref{airy}.
Expanding  $v$ into the Fourier series \eqref{fourier} allows one to formally write 
\begin{equation}\label{solairy}
u(\theta,t):=\sum_{n=-N}^N
    f_n e^{-in^3 t}e^{i n\theta}  \, ,
\end{equation}
which by elementary differentiation of the exponential functions is  seen to solve the Airy equation.

\section{Nonlinear Fourier series}

Throughout the twentieth century and beyond, nonlinear variants of Fourier series have appeared in different contexts. Here, we motivate them within the modern context of representing functions in quantum computing.

\begin{figure}
    \includegraphics[width=0.9\textwidth]{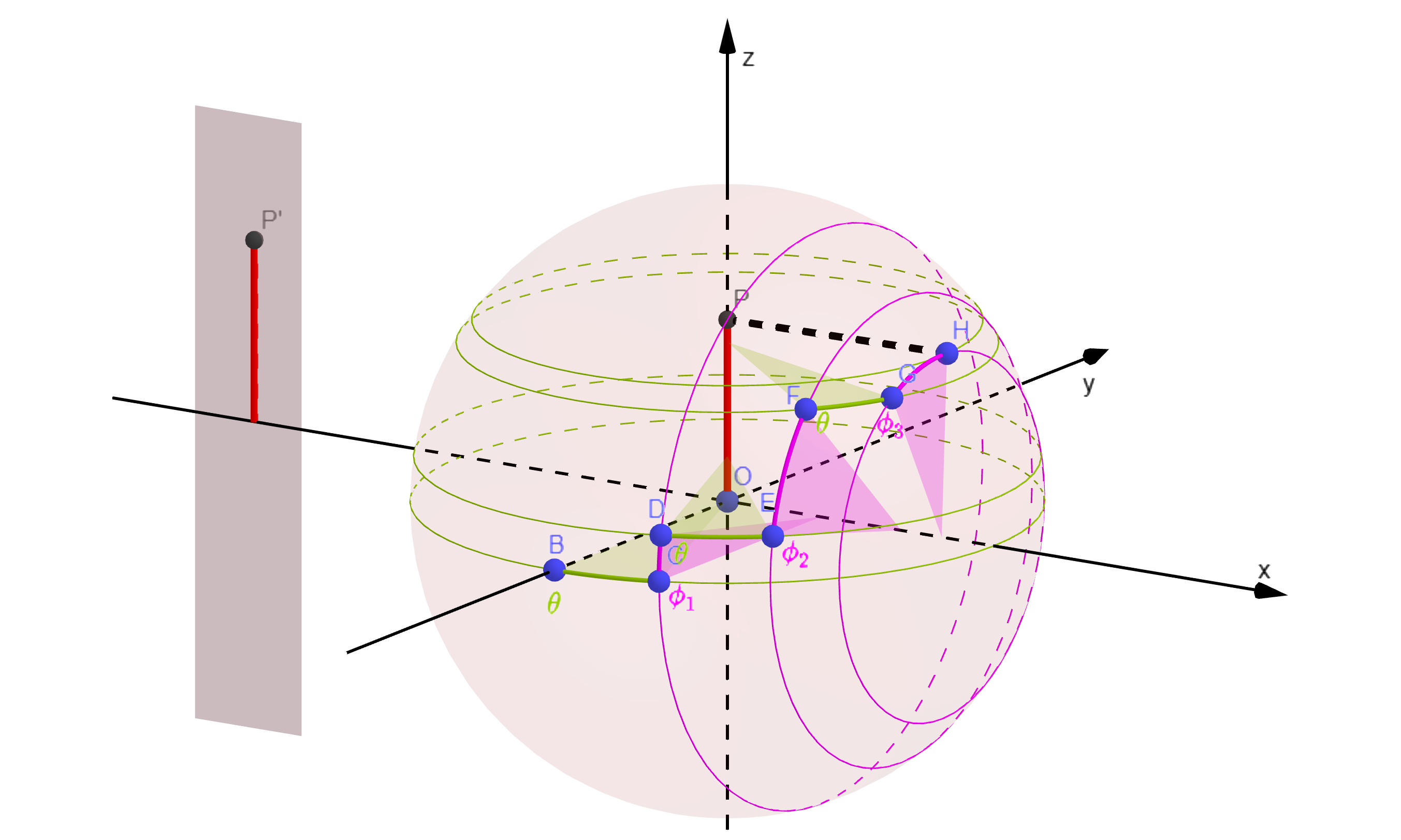}
    \caption{Nonlinear Fourier series}
    \label{NFs}
\end{figure}

The trigonometric polynomials \eqref{fourier} are constructed by composing two basic operations, namely multiplication by $z$ and addition of a complex number. To emphasize the point, we reorganize \eqref{fourier} with the distributive law into an alternation of these two
operations. More precisely, for the Fourier series between
powers $M$ and $N$ we have the following alternation except for a final multiplication by $z^{M}$, 
\begin{equation}\label{distributive}
        \sum_{n={M}}^Nf_n z^n
=        \left[(((f_{N} z +f_{N-1})z+f_{N-2}) z+ f_{N-3})z + \dots +f_{M}\right]z^{M}\, .
    \end{equation}
Multiplication by the complex number $z\in \T$ is a unitary operation on the  Hilbert space $\C$ and thus natural for a quantum computer. However, adding a coefficient $f_n$ is not unitary. For this reason, one deviates from the classical Fourier series and alternates instead
the following two unitary operations in $\C^2$:
\begin{equation}\label{Z}
Z:=\left(\begin{array}{cc} z^{1/2} & 0 \\
0 & \overline{z^{1/2}}\end{array}\right)
\, ,
\end{equation}
\begin{equation}\label{Fn}
F_n:={\frac 1{\sqrt{1+|f_n|^2}}\left(\begin{array}{cc} 1 & f_n \\
-\overline{f_n}  & 1\end{array}\right)}
\, .
\end{equation}
In place of \eqref{distributive},
one considers the following definition of a nonlinear Fourier series with values in the multiplicative matrix group $SU(2)$ for coefficients with indices between $M$ and $N$:
\begin{equation}\label{alternating}
\left(\begin{array}{cc} a_{M,N} & b_{M,N} \\
-\overline{b_{M,N}}  & \overline{a_{M,N}}\end{array}\right):=Z^{{M}}\left(F_MZF_{M+1}Z\dots ZF_{N-1}ZF_N 
\right)Z^{{-N}}\, .
\end{equation}
The alternation of unitary operations is illustrated in  Figure \ref{NFs}. For simplicity, the figure displays the simpler rotations about orthogonal axes in $\R^3$ instead of the $SU(2)$ actions.
The extra left multiplication by $Z^M$ in \eqref{alternating} compared to \eqref{distributive} has the effect that the half-integer powers in \eqref{Z}
multiply to integral powers of $z$ in \eqref{alternating}. For practical computations, this renders obsolete the question of the complex branch of the square root, a point further illustrated by writing \eqref{alternating} in a form
more reminiscent of \eqref{fourier} as
\begin{equation}\label{conjugation}
=
(Z^{{M}}F_MZ^{-M})(Z^{{M+1}}F_{M+1}Z^{-M-1})\dots (Z^{{N}}F_NZ^{-N})
\end{equation}
with non-commuting factors built from iterated conjugation:
\begin{equation}\label{basic}
Z^{{n}}F_nZ^{-n}=\frac 1{\sqrt{1+|f_n|^2}}\left(\begin{array}{cc} 1 & f_nz^n \\
-\overline{f_n}z^{-n} & 1\end{array}\right)\, .
\end{equation}

We observe the following properties of the matrix entries on the left-hand side of \eqref{alternating}.
\begin{theorem}\label{domain}
With the identification $\overline{z}=z^{-1}$ on $\T$, the functions 
\begin{equation}\label{listentries}
    a_{M,N}\,, \ b_{M,N}\, , \ \overline{a_{M,N}}\,, \ \overline{b_{M,N}}
\end{equation}  are Laurent polynomials in $z$
with nonzero coefficients contained in the  closed intervals 
\begin{equation}\label{domainintervals}
    [-N+M,0]\, ,\  [M,N]\, , \ [0,M-N]\, , \ [-N,-M]
\end{equation}
respectively.
\end{theorem}

\begin{proof}
The alternating product within brackets in \eqref{alternating} gives a
matrix-valued Laurent polynomial in $z^{1/2}$ with nonzero coefficients in $[-N+M,N-M]$. Multiplying
the factors outside the brackets shows that the functions \eqref{listentries}
are Laurent polynomials in $z^{1/2}$ with nonzero coefficients in the intervals 
\eqref{domainintervals} multiplied by a factor two.
The expressions \eqref{conjugation} and \eqref{basic} show that 
  \eqref{alternating} is a Laurent polynomial in $z$. As Laurent polynomials in $z$, the functions \eqref{listentries}
  then have coefficients in the intervals \eqref{domainintervals}.
\end{proof}
\begin{theorem}
In first order approximation in $f_n$, $a_{M,N}$ is a constant function equal to one and $b_{M,N}$ is the linear Fourier series 
\begin{equation}\label{MNseries}
\sum_{n=M}^N f_nz^n    \, .
\end{equation}
\end{theorem}
\begin{proof}
    We may ignore the scalar factor in \eqref{basic}, because it is equal to one up to an error of second order in $f_n$. The term \eqref{basic} becomes the identity plus an off-diagonal perturbation 
    of first order in $f_n$. The linear approximation of the product \eqref{conjugation} is then the identity plus the sum of the individual perturbations.
\end{proof}

The limit 
\begin{equation}\label{limitalternating}
\left(\begin{array}{cc} a & b \\
-\overline{b}  & \overline{a}\end{array}\right)=\lim_{N\to \infty}
\left(\begin{array}{cc} a_{-N,N} & b_{-N,N} \\
-\overline{b_{-N,N}}  & \overline{a_{-N,N}}\end{array}\right)
\end{equation}
exists at every point $z$ if the sequence $f_n$ is absolutely summable. This follows from
elementary  convergence criteria for infinite products. If the sequence $f_n$ is merely square integrable, then the limit does
not necessarily exist for every $z$ in the unit circle $\T$. However, in analogy to the linear case, 
the limit exists in the sense of a suitable metric on the space of functions on $\T$ with values in $SU(2)$.

As a surprising difference between the linear and nonlinear cases, the nonlinear Fourier series is not injective on the set of square integrable sequences. In a second surprise, one can compute a unique solution for a related question of interest in quantum signal processing. Here, one is interested in a function $b$,  typically assumed to be real-valued, and wishes to compute a sequence $f_n$ representing $b$ as the upper 
off diagonal entry of the nonlinear Fourier series. The function $a$ has an auxiliary character and can be freely chosen. Choosing $a$ as what is called an outer function, see \cite{MR4974718} for the precise definition, leads to the following theorem.

\begin{theorem}\label{qsp}
For each measurable $v:[0,2\pi)\to [-1,1]$ such that $\log(1-|v(\theta)|^2)$ is integrable, setting $b(e^{i\theta}):=v(\theta)$, there is a unique outer function $a$ and a unique square integrable sequence $f_n$ such that the first row of the nonlinear Fourier series \eqref{limitalternating} of $f_n$
consists of the said $a$ and $b$.
\end{theorem}

Much earlier than any interest in quantum signal processing, a different model of nonlinear Fourier 
series gained prominence in the theory of orthogonal polynomials, see \cite{MR2105088}. Replacing \eqref{Fn} by
\begin{equation}\label{SU(2)}
{F}_n:={\frac 1{\sqrt{1-|f_n|^2}}\left(\begin{array}{cc} 1 & f_n \\
\overline{f_n}  & 1\end{array}\right)}
\end{equation}
with $|f_n|<1$, the correspondingly modified product \eqref{alternating}, where on the left-hand side we also have to drop the minus sign in the lower left entry, produces a nonlinear Fourier series valued in the Lie group
$SU(1,1)$. In the theory of orthogonal polynomials,
we assume that the coefficients $f_n$ are supported on the set of positive integers and set $a_N=a_{0,N}$, $b_N=b_{0,N}$ for $N\ge 0$ using $f_0=0$.
Thanks to the analog of Theorem \ref{domain}, the function
\begin{equation}\label{orthpol}
    \phi_M(z):=z^M(a_M(z)+\overline{b_M(z)})
\end{equation}
is a polynomial in $z$ on $\T$ of degree at most $M$. Extending the product formula $\eqref{conjugation}$
from $0$ to $N$  with $N>M$ and then dividing from the right the 
product from $M+1$ to $N$ expresses \eqref{orthpol} as:
\begin{equation}\label{longop}
    \phi_M(z):=z^M((a_N(z)+\overline{b_N(z)})\overline{a_{M+1,N}(z)}-(b_N(z)+\overline{a_N(z)})\overline{b_{M+1,N}(z)})\, .
\end{equation}
Two polynomials  $\phi_M$, $\phi_{M'}$ with different indices $M,M'<N$ turn out to be orthogonal  with respect to the 
inner product 
\begin{equation}\label{innerprod_N}
   \left< \phi,\psi\right>_N:=   \int_{\mathbb{T}}  \frac{\phi(z)\overline{\psi(z)}}
    {(a_N(z)+\overline{b_N(z)}(\overline{a_N(z)} + b_N(z))} \frac{|dz|}{2\pi}\, .
\end{equation}
This follows by inserting \eqref{longop}
for the polynomial of larger degree 
and identifying with the analog of Theorem \ref{domain} a vanishing zeroth coefficient of the resulting Laurent series in the integrand.
Taking a weak limit towards a positive measure $\mu$,
\begin{equation}\label{innerprod}
   \left< \phi,\psi\right>:= \lim_{N\to \infty} \left< \phi,\psi\right>_N 
    =:\int_{\mathbb{T}}  {\phi(z)\overline{\psi(z)}}\, d\mu(z)\, ,
\end{equation}
gives an inner product so that
all $\phi_M$ are mutually orthogonal.

It turns out that there are analogous polynomials in the $SU(2)$ model, albeit only half orthogonal  \cite{alexis2025sidedorthogonalpolynomialspointwise}. Here, the relevant pairing
\begin{equation*}
   \left< \phi,\psi\right>_N:=   \int_{\mathbb{T}}  \frac{\phi(z)\overline{\psi(z)}}
    {(a_N(z)+\overline{b_N({z})})(\overline{a_N({z})} - b_N(z))}\frac{|dz|}{2\pi}
\end{equation*}
does not define an inner product as the weight is in general not real and the pairing lacks symmetry. One is led to two sequences of polynomials,
\begin{equation*}
    \phi_M(z):=z^M(a_M(z)+\overline{b_M({z})})\, , \  
    {\psi}_M(z):=z^M(a_M(z)-\overline{b_M({z})})\, ,
\end{equation*}
with the relation
\begin{equation*}
    \left< \phi_{M'},{\psi}_M\right> =0
\end{equation*}
for two different $M', M$ smaller than $N$.  Hence $\phi_M$ are left orthogonal polynomials and ${\psi}_M$
are right orthogonal polynomials.

Another application of nonlinear Fourier analysis is to the solution of partial differential equations analogous to Fourier's original idea. 
Following \cite{TaoThiele2012}, we pass to Fourier integrals rather than
series. The integral version of \eqref{fourier} is
\begin{equation}\label{Fourier integral}
\widehat{f}(\xi):=\int_{-\infty}^\infty f(t)e^{-it\xi} dt\, .
\end{equation}
The map from $f$ to $\widehat{f}$ formally coincides with its inverse except for a minus sign in the argument.
We will apply Fourier's idea to the inverse of the Fourier transform.

To motivate the nonlinear variants of \eqref{Fourier integral}, we take a multiplicative view 
 and exponentiate \eqref{Fourier integral}. Define
\begin{equation*}
g(\xi,x):=\exp( \int_{-\infty}^x f(t)e^{-it\xi} dt)\, .
\end{equation*}
The function $g$ formally satisfies for fixed $\xi$ the differential equation
\begin{equation}\label{diffg}
\partial_x g(\xi,x) = g(\xi,x) f(x)e^{- i x\xi}\, 
\end{equation}
with the limiting boundary data 
\begin{equation}\label{boundg}
g(\xi,-\infty)=1\, , \ g(\xi,\infty)=e^{\widehat{f}(\xi)}\, .
\end{equation}
The $SU(2)$ and $SU(1,1)$ Fourier transforms are defined by 
the following matrix-valued analogs of \eqref{diffg} and \eqref{boundg}
with the upper and lower choice of signs, respectively:
\begin{equation}\label{SUcont} \partial_x G(\xi,x)
= G(\xi,x) \left(\begin{array}{cc} 0 & f(x) e^{- ix\xi} \\   \mp \overline{f}(x) e^{ ix\xi} & 0 \end{array}\right) \, ,\end{equation}
\begin{equation}\label{nlft}
    G(\xi,-\infty) =
\left(\begin{array}{cc} 1 & 0 \\
0 & 1\end{array}\right)\, , \ 
   G(\xi,\infty)
=:
\left(\begin{array}{cc} a(\xi) & b(\xi) \\
 {\mp} \overline{b}(\xi) & \overline{a}(\xi)\end{array}\right)\, .
\end{equation}
In linear approximation, the function $a$ is constant
$1$ and the function $b$ is the Fourier transform of $f$.
These nonlinear perturbations of the Fourier transform can be used to solve two nonlinear
perturbations of the Airy equation, namely the  focusing and the defocusing modified Korteweg-de Vries (mKdV) equations 
\begin{equation}
    \partial_t u(x,t)= \partial_x^3 u(x,t)   \pm 6u^2(x,t)\partial_xu(x,t)
\end{equation}
with the initial condition 
\begin{equation}
    u(x,0)= f(x)\, .
\end{equation}
If $a(\xi)$ and $b(\xi)$ denote the entries in the first row of the $SU(2)$ Fourier transform of $f$ as in \eqref{nlft},
then $a(\xi)$ and $b(\xi)e^{i\xi^3t}$ are the entries in the first row of the $SU(2)$ Fourier transform of the solution $u(x,t)$ of the 
focusing mKdV equation at time $t$. Similarly, we use the $SU(1,1)$ Fourier transform to determine the solution to the defocusing equation.
The multiplication of the entry $b$ by a phase $e^{i\xi^3t}$ is analogous to the solution to the Airy equation \eqref{airy}
presented in the introduction, with the roles of $f$ and $\widehat{f}$ interchanged.
A rigorous justification of this solution method requires a bit of work and in particular an analog of an existence
and uniqueness result for the inverse nonlinear Fourier transform under suitable assumptions. 

The mKdV equations are named after two Dutch mathematicians, the first of whom, Korteweg, was the PhD advisor both
of de Vries and Brouwer. The closely related original KdV equation was the first to explain solitary waves as
observed in shallow water channels
\cite{MR3363408}.

\section{Fast Fourier transforms}
The numerical evaluation of the Fourier series \eqref{fourier} is limited to finitely many points $z$ of $\T$. 
The natural choice are equidistant points, say $2^m$ many. We always assume that $1$ is one of these equidistant points. The  values  at these sampling points do not distinguish between two exponentials of \eqref{fourier} whose indices differ by a multiple of $2^m$. 
We therefore restrict the sum of \eqref{fourier} to the half open interval $I:=[0,2^m)$.
The celebrated fast Fourier transform (FFT) evaluates \eqref{fourier} at the $2^m$ sampling points in
$O(m2^m)$ computational steps. By analogy between the forward and inverse problems for finite Fourier series, the same algorithm can 
be used to compute the coefficients when the values at the sampling points are given.

In the nonlinear setting, the forward problem of evaluating a nonlinear Fourier series \eqref{alternating} at the sampling points 
for given coefficients $f_n$ is different from the inverse problem of computing the coefficients
given the values of \eqref{alternating} at the sampling points. We present algorithms for both tasks, each requiring $O(m^22^m)$ computational steps.
 
\subsubsection*{The forward algorithm} Let coefficients $f_n$ with $n\in I$ be given. A subinterval $J$ of $I$ is called dyadic, if it is  
of the form \begin{equation}
    J=[2^kn,2^k(n+1))
\end{equation} 
with nonnegative integers $k$ and $n$. Dyadic intervals may be arranged in a tree based on inclusion, see Figure \ref{example_tree}.
The algorithm computes for every such $J$ at  $2^k$ equidistant points
the nonlinear Fourier series  \eqref{alternating} with $M$ and $N$ given by 
\begin{equation*}
    [M,N+1)=J\, .
\end{equation*}
 We call this nonlinear Fourier series $\widehat{J}$.
 In particular, the algorithm computes $\widehat{I}$ at the desired sampling points.

The algorithm sorts the dyadic subintervals $J$
by increasing $k$.
If $k$ is zero, then \eqref{conjugation} has one factor and we read from \eqref{Fn} the matrix entries of $\widehat{J}$ sampled at $z=1$.
If $k>0$, then the left half $L$ and the right half $R$ of $J$ are also dyadic. They are called the children of $J$. As the children are shorter than $I$, the algorithm has already computed  $\widehat{L}$ and $\widehat{R}$. The product representation \eqref{conjugation} gives
\begin{equation}\label{FT multiplication}
    \widehat{J}=\widehat{L}\widehat{R} \, .
\end{equation}
More precisely, the left-hand side of \eqref{FT multiplication} requires twice the number of sampling points as computed
for $\widehat{L}$ and $\widehat{R}$. Therefore, we expand the latter into linear Fourier series while
resolving the periodic ambiguity of the linear Fourier modes by  placing the modes of
each matrix entry of $\widehat{L}$ and $\widehat{R}$ into
their natural domains, by which we mean the intervals described in Theorem \ref{domain} for a given interval.
Then we evaluate these linear Fourier series at the sampling points required for $\widehat{J}$ and compute the product \eqref{FT multiplication}. 
This procedure requires $O(k2^k)$ computations  using the FFT. Summing over the
$2^{m-k}$ intervals at scale $k$ and then summing over the $m$ scales bounds the total number of 
computations by $O(m^22^m)$.

\subsubsection*{The inverse algorithm}

It is known, 
that the coefficients $f_n$ in the interval $I$ are uniquely determined by the values of 
$\widehat{I}$ at $2^m$ equidistant sampling points \cite{tsai}.  
We present an algorithm that computes the coefficients $f_n$ following \cite{ni2025inversenonlinearfastfourier}, which was inspired by \cite{MR1011737}. 
Like the forward algorithm, this inverse algorithm computes $\widehat{J}$ at $2^k$ equidistant sampling points for each dyadic subinterval $J$ of $I$.
In addition, it computes an auxiliary matrix-valued function $\overline{J}$ at the same $2^k$ sampling points. 

\begin{figure}
\centering
  \includegraphics[width=0.5\linewidth]{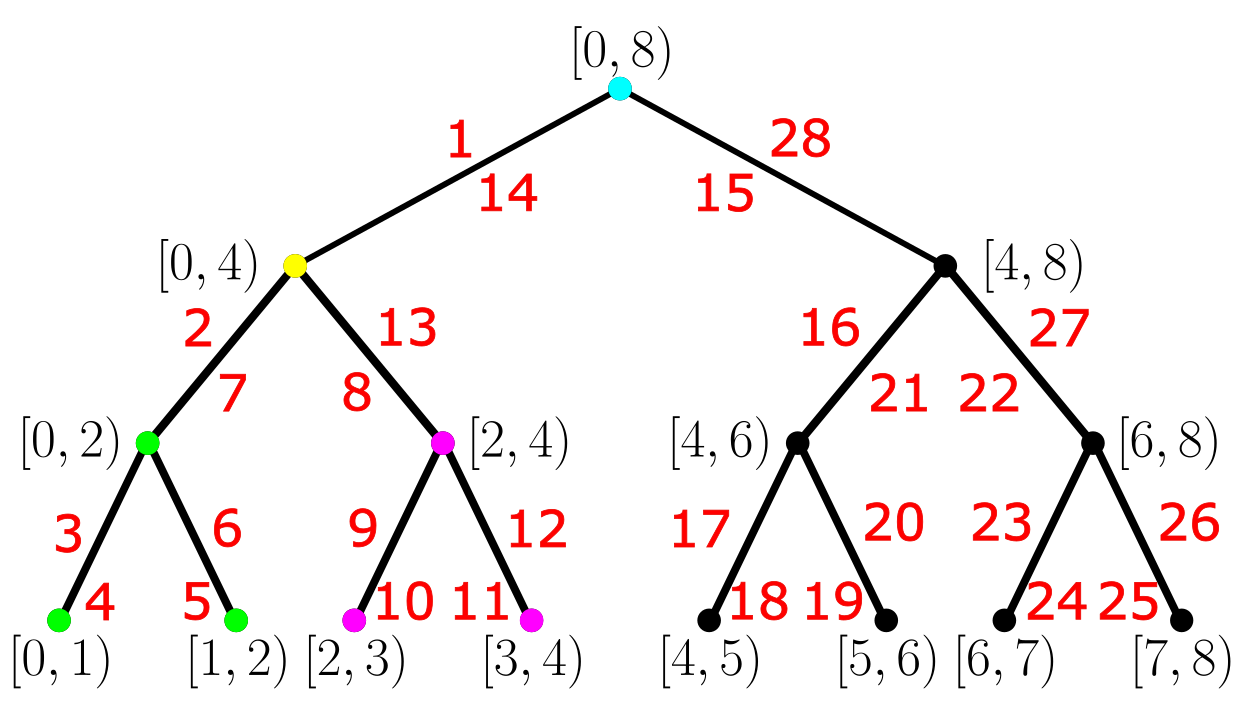}
  \caption{The tree of dyadic subintervals of $[0, 8)$. The edges are enumerated in the order they are visited. 
The yellow interval $J$ is visited three times. 
Prior to each visit to $J$, the algorithm has visited all the blue intervals, all the blue or green intervals, and all the blue or green or magenta intervals, respectively.}
  \label{example_tree}
\end{figure}

The algorithm runs along a path through the tree of dyadic intervals, shown in Figure \ref{example_tree} in case $m=3$.
Each node of the tree other than those in the bottom layer is visited three times, and each edge is once traveled down and later traveled back up, in the order shown in the figure. 

When the algorithm reaches a dyadic interval $J$ with children for the first time,
the algorithm computes $\overline{L}$ for the left child $L$ of $J$ and then proceeds to $L$.
When it reaches the interval $J$ for the second time, it computes $\overline{R}$ for the right child $R$ of $J$
and then proceeds to $R$. When it reaches the
interval $J$ for the third time, it computes $\widehat{J}$ and returns to the parent of $J$.
If the algorithm reaches an interval $J$ without children, then the algorithm  computes directly $\widehat{J}$ and returns to the parent. The algorithm starts and ends with the 
full interval $I$.  

In more detail, we initialize $\overline{I}$ to be equal to the given $\widehat{I}$.
Assume that we have reached an interval $J$ for the first time and that $J$ has children. Unless $J=I$, we arrived from the parent of $J$. In either case, we know $\overline{J}$ and compute
\begin{equation}\label{project left}
\overline{L}:=P_L\overline{J}\, .
\end{equation}
Here, the operation $P_L$ is a linear Fourier multiplier that multiplies the linear 
Fourier coefficients of the matrix entries of $\overline{J}$, assumed to be supported in the natural domains  associated with $\widehat{J}$, by the indicator function of the natural domains associated with $\widehat{L}$. We  then sample these so restricted linear Fourier series at the sampling points associated with $\overline{L}$.

When we reach $J$ for the second time, we have returned from the left child and
we know $\widehat{L}$ in addition to $\overline{J}$. We then compute
\begin{equation}\label{polynomial cutoff}
    \overline{R}:=P_R (\widehat{L}^{-1}\overline{J})\, .
\end{equation}
More precisely, we first compute the linear Fourier series of $\widehat{L}$ and $\overline{J}$ in the natural domains. Then we compute the
linear Fourier series of $\widehat{L}^{-1}$ by reshuffling the entries of $\widehat{L}$, using that the determinant of the values
of $\widehat{L}$ are one. Then we evaluate the Fourier series at the sampling points associated with $\overline{J}$ and
take the product \eqref{polynomial cutoff}. We expand this product in linear Fourier series with modes in the natural
domains of $\overline{J}$. Then we apply the Fourier multiplier $P_R$, defined analogously to $P_L$ above, and sample the restricted Fourier series at the sampling points associated with $\overline{R}$.

When we reach the interval $J$ for the third time, we have returned from the right child and we know 
$\widehat{R}$ and $\widehat{L}$. We compute $\widehat{J}$ as in \eqref{FT multiplication}.

When we reach an interval $J$ that contains only one integer $n$, we compute $f_n:=b_n \overline{a_0}^{-1}$, where $a_0$ and $b_nz^n$ are the entries in the first row of $\overline{J}$.
Then we compute $\widehat{J}$ with \eqref{Fn}.

This completes the description of the algorithm. Note that the coefficients
$f_n$ are obtained from the computations at the intervals of length one, and we may also stop the algorithm after computing all these.  The algorithm  relies on a process called `layer stripping'.  This first factor in the product \eqref{conjugation} is uniquely determined by the property
that dividing the product from the left by this factor results in a shorter product. This property leads to the expression $f_M=b_M\overline{a_0}^{-1}$, where $b_M$ is the $M$-th coefficient of $b_{M,N}$ and $a_0$ is the constant coefficient of $a_{M,N}$. If one replaces  the projections $P_L$ and $P_R$ in the above algorithm by identity operations, the first nontrivial step of the algorithm is to determine the first factor of $\eqref{conjugation}$
by the layer stripping identity. The next step is to divide $\eqref{conjugation}$ from the left by this factor. The subsequent steps iterate this procedure with the remaining factors.

One step of layer stripping only needs a single linear Fourier mode of each of $a_{M,N}$ and $b_{M,N}$. Similarly, a number of iterative steps of layer stripping only need as many linear Fourier modes of
$a_{M,N}$ and $b_{M,N}$ as there are steps. The projections $P_L$ and $P_R$ in the algorithm project at every step to
the necessary Fourier modes. This reduces the number of computational steps to $O(m^22^m)$, as one sees by counting
similarly to the forward algorithm.

\section{Acknowledgements}

The second author was supported by the ERC Synergy Grant 101224275.
\printbibliography
\end{document}